\documentclass{rspublic}
\pdfoutput=1
\usepackage[numbers]{natbib}
\usepackage[font=footnotesize]{caption}
\usepackage{url,color,graphicx,amsmath,amsfonts}
\usepackage[color]{changebar}
\cbcolor{red}
\usepackage[normalem]{ulem}
\usepackage{mathptmx}

\usepackage{microtype}

\newcommand{\R}{\mathbb{R}}
\newcommand{\Q}{\mathbb{Q}}
\newcommand{\id}{I}
\newcommand{\iu}{\mathrm{i}}

\renewcommand{\d}{\mathop{}\!\mathrm{d}}
\renewcommand{\epsilon}{\varepsilon}
\renewcommand{\mod}{\operatorname{mod}}
\newcommand{\im}{\operatorname{Im}}
\renewcommand{\re}{\operatorname{Re}}

\title[Bifurcation analysis of delay-induced resonances of ENSO]
{Bifurcation analysis of delay-induced resonances of the El-Ni{\~n}o
  Southern Oscillation}

\author[B. Krauskopf and J. Sieber]{Bernd Krauskopf$^1$ and Jan Sieber$^2$,
}

\affiliation{
  $^1$Department of Mathematics,
The University of Auckland, New Zealand\\
  $^2$College of Engineering, Mathematics and Physical Sciences,
  University of Exeter, UK}

\begin{document}

\maketitle

\begin{abstract}{delay differential equations, bifurcation analysis,
    El-Ni{\~n}o Southern oscillation} 

    Models of global climate phenomena of low to intermediate
    complexity are very useful for providing an understanding at a
    conceptual level. An important aspect of such models is the
    presence of a number of feedback loops that feature considerable
    delay times, usually due to the time it takes to transport energy
    (for example, in the form of hot/cold air or water) around the
    globe. 

    In this paper we demonstrate how one can perform a bifurcation
    analysis of the behaviour of a periodically-forced system with
      delay in dependence on key parameters. As an example we
    consider the El-Ni{\~n}o Southern Oscillation (ENSO), which is a
    sea surface temperature oscillation on a multi-year scale in the
    basin of the Pacific Ocean. One can think of ENSO as being
    generated by an interplay between two feedback effects, one
    positive and one negative, which act only after some delay that is
    determined by the speed of transport of sea-surface temperature
    anomalies across the Pacific.  We perform here a case study of a
    simple delayed-feedback oscillator model for ENSO (introduced
    by Tziperman \emph{et al}, J. Climate 11 (1998)), which is
      parametrically forced by annual variation.  More specifically,
    we use numerical bifurcation analysis tools to explore directly
    regions of delay-induced resonances and other stability boundaries
    in this delay-differential equation model for ENSO.
\end{abstract}

\section{Introduction}
\label{sec:intro}

In climate science one finds a hierarchy of increasingly complex
mathematical models, all the way from low-dimensional conceptual
models in the form of ordinary differential equations (ODE) to global
climate models with many millions of unknowns. Our interest here is in
models of a low to intermediate complexity, in particular, those that
feature delayed feedback mechanisms. The mathematical description of
such a system, hence, takes the form of delay differential equations
(DDEs). Despite looking similar to ODEs at first sight, DDEs have
an infinite-dimensional phase space (that is, space of possible
initial conditions), which enables them to describe complex
phenomena even with a low number of dependent variables
\cite{AK98}. This type of
model arises in other fields of science, for example, in mechanical
engineering, in laser physics and in mathematical biology, where one
finds communication delays (between optical elements or cells) as well
as processing delays; see, for example, \cite{S89,K05}. In climate
science, on the other hand, delays arise mainly in feedback loops due
to the transport of mass or energy from one location on the planet to
another. There are models that describe such transport phenomena
directly in the form of partial differential equations
\cite{D08}. However, in some situations it is advantageous to take a
more conceptual point of view by considering only the resulting
feedback effect itself, which is then subject to a delay due to the
associated transport time.

In this contribution we perform a case study of such a DDE model to
demonstrate how tools from bifurcation theory can be brought to bear
to understand its behaviour. More specifically, we study delay-induced
resonances in the classical ENSO problem. This is an ideal showcase
problem for bifurcation analysis because the ocean-atmosphere system
of the Pacific has an oscillatory instability at the linear level due
to reasonably well identified feedback effects.  The resulting
oscillator interacts with the annual forcing to produce a rich set of
resonances \cite{JNG94,JNG96,TSCJ94} in a wide region of parameters.
Extensive research has created a hierarchy of models that isolate the
processes affecting ENSO;
see also \cite{NBx98}. We explore here one of the simplest of these
models that still contains physically meaningful parameters, namely,
the model introduced by Tziperman \emph{et al} \cite{TCZXB98}.  This
model takes the form of a scalar DDE that lumps all coupled processes
into two feedback terms, which act only after given delays due to the
time it takes for thermocline anomalies (proportional to sea-surface
temperature anomalies) to travel across the Pacific.  We provide a
brief summary of the motivation behind this delayed-feedback oscillator model
of ENSO in Section~\ref{sec:model}; for an account of the modelling of
the processes involved at coarse to intermediate levels of complexity
see the review \cite{NBx98} and the textbook \cite{D08}.

Once a DDE model, such as that for ENSO, has been derived, the task is
to understand its behaviour in dependence on relevant parameters.  A
quite straightforward method for investigating the eventual or
attracting behaviour of a given mathematical model is its simulation
by numerical integration. However, when one wants to find all possible
behaviour by simulation, one is faced with a number of challenges.
One generally needs to perform many simulation runs for different
parameter values, for example, chosen from some sufficiently fine
grid; each such run needs to be performed until transients die out
sufficiently and care needs to be taken in the choice of the initial
condition, especially in the presence of multistability. Note that for
DDEs one needs an entire history segment as initial condition
\cite{HL93}, making the search space infinite dimensional. Furthermore, there
is the difficulty of representing and identifying, ideally
automatically, what the eventual dynamics actually is.

An alternative to simulation is \textit{bifurcation analysis}.  It was
orginally developed as a mathematical tool to classify the long-term
behaviour of nonlinear dynamical systems that are modelled by
low-dimensional ordinary differential equations (ODEs) or maps. The
underlying idea is that one studies equilibria and periodic orbits and
their stability and other properties as a function of relevant system
parameters. When a parameter is varied, one may observe well-defined
qualitative changes of the dynamics, which are referred to as
bifurcations. The goal of bifurcation analysis is to find the
\emph{bifurcation diagram}, which divides the parameter space into
regions of qualitatively different behaviour. A bifurcation diagram is
effectively a `road map' of where in parameter space one may find
which kind of behaviour, and what changes occur in the transition from
one type of behaviour to another. In spite of their seeming
simplicity, ODE models of quite low dimension (such as the well-known
Lorenz and R\"ossler systems) can already give rise to intricate
bifurcation diagrams with regions of stationary, periodic,
quasiperiodic and even chaotic behaviour; see, for example, the
textbooks \cite{TS02,GH83,K04} as entry points to bifurcation theory.

Today, the bifurcation analysis of a given ODE can be performed
routinely with freely available numerical tools, such as the packages
\texttt{AUTO} \cite{DCFKSW98} and \texttt{MatCont} \cite{DGK03}; see
also the survey \cite{GK07}.  These tools only require knowledge of
the right-hand side of the ODE and allow the user to perform a
systematic exploration of the bifurcation diagram. This is achieved by
finding and then tracking (or continuing) equilibria, periodic orbits,
and their bifurcations to yield stability boundaries of the different
types of solutions.  The usefulness of bifurcation theory has been
demonstrated in numerous studies of ODE models arising in diverse
fields of science; see, for example, \cite{KOG07}.

More recently, numercial bifurcation tools have been implemented also
for more general classes of systems. For complex models, which are
typically based on partial differential equation formulations,
numerical routines are available as libraries
\cite{LOCA02}. Review \cite{Da14} discusses bifurcation analysis
techniques specifically for fluid dynamics problems (which are the
dynamical core of climate models). Textbook \cite{D08} gives an
introduction to bifurcation analysis in oceanography and
demonstrates the role it can play (it uses ENSO as one of its case
studies). Specifically for DDEs, numerical methods for their
bifurcation analysis are today available in the form of the packages
\texttt{DDE-Biftool} \cite{ELR02} and \texttt{knut} \cite{SSH06}.
These numerical tools apply to general systems of DDEs, including DDEs
with several delays and periodic forcing.  In the same way as for
ODEs, these tools can determine and track the stability properties and
bifurcations of equilibria and periodic orbits.  As we will
demonstrate, this means that numerical bifurcation studies can be
performed efficiently today also for models that feature delays; see
also the survey \cite{RS07}.

The most obvious application areas of bifurcation analysis are those
where systems under consideration have tunable parameters that the
practitioner can adjust to achieve the desired behaviour \cite{TS02};
concrete such examples are engineering and laser physics.  The
situation in climate science, on the other hand, is very different in
that tuning a parameter is typically not possible in the real
system. In spite of this fact, we argue that systematic parameter
studies are still relevant in this context, namely to determine the
sensitivity of the model under consideration to changes of parameters
that are known only in an order-of-magnitude sense. The fact that the
observed behaviour of a nonlinear system may depend critically on all
parameters motivates the bifurcation analysis with respect to selected
difficult-to-determine parameters; see \cite{DG05} for an exemplary
study.

For our specific example, two parameters that affect the ENSO
resonances strongly are the mean and the annual variation of the
ocean-atmosphere coupling. The values of these two parameters are hard
to ascertain in the real system; furthermore, they are also difficult
to compare and convert between models of different complexity. Hence,
it makes perfect sense to locate a chosen resonance tongue in the
given model to restrict the values of these two parameters. When one
attempts to locate a particular resonance one encounters the problem
that the only places where one knows the frequency ratio analytically
(or at least numerically without a scan of the complete parameter
plane) are the small-forcing and the small-oscillation-amplitude
regimes. In these two regimes resonance tongues are extremely narrow,
and convergence to the resonant orbits (if they are stable) is very
slow. This makes resonances hard to observe and track with simulations
in the small-amplitude regime (corresponding to the classical devil's
staircase scenario \cite{JNG94}).  In this paper we show how this
difficulty can be overcome with numerical bifurcation analysis. More
specifically, periodic orbits are tracked as solutions of a
boundary-value problem (BVP) and, as such, these computations are
unaffected by the instability (or weak stability) of the periodic
orbit under consideration or by its sensitive dependence on
parameters. With this BVP approach we have tracked the relevant
resonance tongues directly from their root point near the linear
regime over a considerable region of the parameter plane.  This
bifurcation diagram, shown in Figure~\ref{fig:bif2d} of
Section~\ref{sec:bifs} below, constitutes our main result which
demonstrates the capabilities of numerical bifurcation analysis in
the context of climate science.

\section{ENSO as a parametrically
  forced delayed-feedback oscillator}
\label{sec:model}

We consider a conceptual model for the ENSO mechanism that simplifies
the circulation of ocean surface water waves and the ocean-atmosphere
interaction to a scalar DDE, which is periodically forced due to the
annual variation of the strength of ocean-atmosphere coupling. We use
the parameters and follow the notation of \cite{TCZXB98}, where the
dependent variable is the height anomaly $h(t)$ of the thermocline at
the Eastern boundary of the Pacific. The height anomaly $h$ (in
dimensionless units) satisfies
\begin{align}
  \label{eq:tzip}
  \frac{\d}{\d\, t} h(t)&=b G\left[\kappa(t-\tau_1)h(t-\tau_1)\right]-
  c G\left[\kappa(t-\tau_2)h(t-\tau_2)\right]-d h(t)\mbox{, \ where}\\
  G(x)&=
  \begin{cases}
    b_+\tanh(x/b_+) & \mbox{if $x\geq0$,}\\
    b_-\tanh(x/b_-) & \mbox{if $x<0$,}
  \end{cases}\label{eq:gdef}\\
  \kappa(t)&=k_0+d_k\sin(\omega t)\label{eq:kappadef}\mbox{.}
\end{align}
The right-hand side of \eqref{eq:tzip} combines positive feedback due
to the warm equatorial sea-surface temperature (SST) perturbation in
the central Pacific, negative feedback due to the cold off-equatorial
SST perturbation in the central Pacific, and dissipation. The
parameters and their values used here are from Tziperman \emph{et
al} \cite{TCZXB98}, and they are listed in Table~\ref{tab:parameters}.

\begin{table}
  \centering
  \begin{tabular}[t]{clp{0.45\textwidth}}\hline\noalign{\smallskip}
    parameter & value (units) & interpretation\\\noalign{\smallskip}\hline\noalign{\smallskip}
    $\tau_1$ & 1.15 (months) &  time it takes warm equatorial SST perturbations 
    from central Pacific to reach eastern boundary  (via Kelvin wave), \\
    $\tau_2$ & 5.75 (months) & time it takes cold off-equatorial SST peturbations from 
    central Pacific to reach eastern boundary  (first westward via
    Rossby wave, then eastward via reflected Kelvin wave), \\
    $b$ & $1/120$ (per day, $0.2535$ per month) &  amplification factor of warm perturbations,\\
    $c$ & $1/160$ (per day, $0.1901$ per month) &  attenuation factor of cold perturbations,\\
    $d$ & $1/190$ (per day, $0.1601$ per month) & attenuation  due to dissipation,\\
    $G(\kappa h)$ & defined in Eq.~\eqref{eq:gdef} & ocean-atmosphere coupling (saturation nonlinearity),\\
    $\kappa(t)$ & $\kappa(t)=k_0+d_k\sin(\omega t)$ & annually varying slope of $G$ at $h=0$,\\
    $\omega$ & $2\pi/12$ (per month) & frequency of annual variation (period  $T=12$ months),\\
    $b_\pm$ & $b_+=3$, $b_-=-1$ & maximum and minimum of saturation nonlinearity in $G$.
    \\\noalign{\smallskip}\hline
  \end{tabular}
  \caption{Values and interpretation of parameters in \eqref{eq:tzip}--\eqref{eq:kappadef}}
  \label{tab:parameters}
\end{table}

\begin{figure}[t!]
  \centering
  \includegraphics[width=1\textwidth]{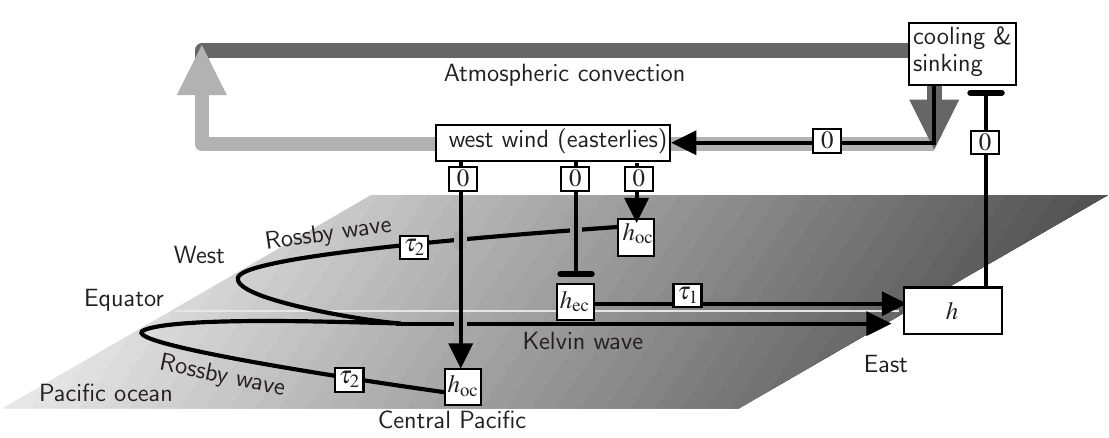}
  \caption{Schematic of the two feedback mechanisms that motivate the
    ENSO model \eqref{eq:tzip}. The shading of the ocean corresponds
    to the thermocline (and sea-surface temperature) in equilibrium:
    darker means lower (colder). Equation~\eqref{eq:tzip} only
    describes the anomaly of the thermocline, $h$, at the eastern
    boundary. The labels $h_\mathrm{ec}$ and $h_\mathrm{oc}$ refer to
    the thermocline anomaly in the equatorial and the off-equatorial
    central pacific, respectively, which feed into the feedback
    mechanisms.  Regular arrows indicate same-sign re-inforcement and
    horizontal (inhibition) bars indicate oppposite-sign
    re-inforcement; the label on each arrow shows the time delay of
    the interaction.}
  \label{fig:ensomech}
\end{figure}

The idea that a dynamical instability is the basic mechanism behind
ENSO was first proposed by Bjerknes in 1969 \cite{B69}. The review by
Neelin \emph{et al} \cite{NBx98} describes the hierarchy of models
that have been developed during the following decades. The
  textbook \cite{D08} presents a detailed walk-through of this hierarchy,
  including a derivation of delayed-feedback oscillator models very similar to
  \eqref{eq:tzip}. The more complex models in this hierarchy help to
establish which physical processes are connected. The simpler models,
such as \eqref{eq:tzip} considered here, take these connections for
granted and treat them as lumped positive or negative feedback terms
in the right-hand side. Figure~\ref{fig:ensomech} sketches the
physical processes that enter the right-hand side of \eqref{eq:tzip}
and their connections.  The height of the thermocline is a proxy for
the sea-surface temperature. In equilibrium the thermocline is low at
the eastern boundary of the Pacific and high at the western end. This
difference maintains an atmospheric convective loop above the
equatorial ocean, which points westward at the ocean surface because
the relatively cold water at the eastern end cools the atmosphere,
causing air to sink (and vice-versa, the relatively warm water at the
western end heats up the atmosphere, causing air to rise). The dependent
variable $h$ in \eqref{eq:tzip} measures the anomaly of the
thermocline at the eastern boundary, that is, the deviation from
equilibrium at this point.

A disturbance in the atmospheric convective loop, that is, a slowing
down of the westward surface wind (creating an eastward wind stress of the
ocean surface) in the central Pacific causes surface water transport
toward the equator. Surface water is warm, causing an increase of $h$
in the equatorial central Pacific, and a decrease of $h$ in the
off-equatorial central Pacific. The positive thermocline disturbance
at the equator travels eastward along the equator (in a \emph{Kelvin
wave}), and hits the eastern boundary after time $\tau_1$. The
corresponding increase in sea-surface temperature at the eastern
boundary slows down the atmospheric convection (the relatively cold
water at the eastern boundary is a driver of the convection). This
creates the first, positive feedback term in \eqref{eq:tzip} with the
delay $\tau_1$ (note that $b>0$).

The off-equatorial depression of $h$ travels westward and toward the
equator (in a \emph{Rossby wave}), reflects at the western Pacific
boundary, and then travels as a Kelvin wave eastward, reaching the
eastern boundary after time $\tau_2$. The depression of $h$
corresponds to a cooling down of the sea surface such that atmospheric
convection is enhanced, creating the negative delayed feedback term in
\eqref{eq:tzip} with the delay $\tau_2$ (note that $c>0$).

\begin{figure}[t!]
  \centering
  \includegraphics[scale=0.7]{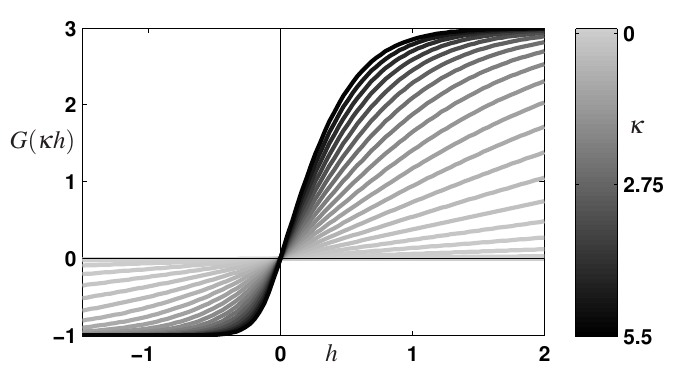}
  \caption{Profile of coupling function $G(\kappa h)$ for $\kappa \in
    0, \dots, 5.5$, which are the values that are encountered in the
    parameter range $(d_k,k_0)\in[0,2.5]\times[0,3]$.}
  \label{fig:gdef}
\end{figure}

The above effects work all at the linear level of the disturbance of
the equilibrium. The nonlinearity in \eqref{eq:tzip} comes from
saturation and asymmetry in the ocean-atmosphere coupling instability,
as expressed by the function $G$ in \eqref{eq:gdef} that is shown in
Figure~\ref{fig:gdef}. Moreover, the slope of $G$ at $h=0$ varies
throughout the year, being strongest in the northern summer and
weakest in the winter; see \cite{TCZXB98} for a realistic time profile
of the annual variation. This annual variation is expressed in the
sinusoidal variation \eqref{eq:kappadef} of the pre-factor $\kappa$ in
the argument of $G$. Since only the amplification factor varies in
time, the oscillator described by \eqref{eq:tzip} is
\emph{parametrically excited}; in particular, the equilibrium $h=0$ of
the unforced system (for $d_k=0$) corresponds to a trivial periodic
solution in the periodically-forced system (for $d_k > 0$).

\section{Delay differential equations and 
  their bifurcations analysis}
\label{sec:genddebif}

This section is a short review of the general theory of dynamical
systems with delays as needed for the bifurcation analysis of the ENSO
model in Section~\ref{sec:bifs}, which the reader may wish to consult
already for illustrations of the concepts introduced now.  To keep the
exposition specific to the ENSO model under consideration, we present
the theory for the case that the right-hand side of the DDE depends on
the current value of the dependent variable, as well as on two values
from the past, with delays $\tau_1\leq\tau_2$. Furthermore, we allow
the right-hand-side of the DDE to depend explicitly on time, which
accounts for the periodic forcing of \eqref{eq:tzip}. Hence, we
consider a DDE of the form
\begin{equation}
\label{eq:gendde}
   \frac{\d}{\d t}x(t) = f(t,x(t),x(t-\tau_1),x(t-\tau_2),\eta),
\end{equation}
where $x \in \R^n$ is an $n$-dimensional state and $\eta\in \R^m$ is a
multi-parameter; one also refers to $\R^n$ as the physical space.
The function 
\begin{equation*}f: \R \times \R^n \times \R^n \times \R^n \times \R^m \to \R^n
\end{equation*}
describes the right-hand-side, and it is assumed to be sufficiently
smooth. Note that \eqref{eq:gendde} is non-autonomous because $f$
depends explicitly on time $t$.  We consider here only the case of
periodic forcing, that is, $f(t,\cdot,\cdot,\cdot,\eta) = f(t +
T_f,\cdot,\cdot,\cdot,\eta)$ where $T_f$ is the forcing period. An
important special case is that the forcing amplitude is zero, so that
system \eqref{eq:gendde} is autonomous. We assume here that the
forcing amplitude is one of the components of the parameter vector
$\eta$, so that the limiting autonomous case can be found as a subset
of parameter space.

We proceed to discuss the basic properties of the DDE
\eqref{eq:gendde}, both for the autonomous and the non-autonomous
cases, which form the basis of the bifurcation analysis of the ENSO
model in Section~\ref{sec:bifs}. We remark that the theory presented
here is also valid for more general DDEs, including those with
distributed delays, as long as all delays are constant and bounded.
Full technical proofs can be found in the classical textbooks by Hale
and Verduyn-Lunel \cite{HL93} and Diekmann \emph{et al} \cite{DGLW95};
the textbook by St{\'e}p{\'a}n \cite{S89} is also rigorous, but
provides a shorter route to those parts of the theory that are
necessary for practical applications.

Because of the presence of delayed terms, the DDE \eqref{eq:gendde}
does not define a dynamical system on the physical space
$\R^n$. Rather, its phase space is the infinite-dimensional space $\R
\times C([-\tau_2,0];\R^n)$; here $\R$ represents time and
$C([-\tau_2,0];\R^n)$ is the space of continuous functions over the
maximal delay interval $[-\tau_\mathrm{max},0] = [-\tau_2,0]$ with
values in the physical space $\R^n$. This means that in order for the
initial-value problem to be well defined one needs to prescribe an
entire function $x_0:[-\tau_2,0]\to\R^n$ as the initial
condition. Similarly, after solving \eqref{eq:gendde} from time $0$ to
time $t$ the current state is the function $x_t:[-\tau_2,0]\to\R^n$
given by $x_t(s)=x(t+s)$,\ $s\in[-\tau_2,0]$.  These trajectories
$x_t$ in the infinite-dimensional phase space $\R\times
C([-\tau_2,0];\R^n)$ depend smoothly on the initial condition $x_0$,
meaning that classical dynamical systems theory can be applied to the
DDE \eqref{eq:gendde} in $\R\times C([-\tau_2,0];\R^n)$. When the
forcing amplitude is zero the time axis decouples and it is sufficient
to study the autonomous DDE with phase space $C([-\tau_2,0];\R^n)$.

\subsection{Equilibria and their bifurcations}
\label{sec:equil}

\noindent
We first consider the case that $f$ is autonomous, which is of
interest on its own as well as in the context of periodic forcing
where it corresponds to zero forcing amplitude. An equilibrium
or steady state is an initial function $x_0:[-\tau_2,0] \to \R^n$ in
phase space such that $x_t(s)=x_0(s)$ for all $s\in[-\tau_2,0]$ and
$t>0$. Hence, $x_0(s)$ equals a constant $\bar x_0$ for all $s\in
[-\tau_2,0]$ and the constant vector $\bar x_0 \in \R^n$ can be determined
from the equation
\begin{equation}
\label{eq:eqdef}
   f(\bar x_0,\bar x_0,\bar x_0,\eta) = 0,
\end{equation}
which, as for ODEs, is a system of $n$ algebraic equations. This means
that equilibria of DDEs can be found and tracked in a parameter as
roots of \eqref{eq:eqdef} by standard continuation methods.

To determine the stability of an equilibrium $x_0$ one moves it to the
origin and linearises the DDE. The result is a linear DDE of the form
\begin{displaymath}
  \dot x(t)= 
  A x(t) +B x(t-\tau_1)+Cx(t-\tau_2)\mbox{,}
\end{displaymath}
where $A$, $B$ and $C$ are $n \times n$ matrices. The stability of the
origin is determined by the spectrum of the eigenvalues of the
linearisation, which can be found as the roots of the
\emph{characteristic function}
\begin{equation}\label{eq:chidef}
  \chi(\lambda)=\det(\lambda\id-A-B e^{-\lambda\tau_1}-
  C e^{-\lambda\tau_2})\mbox{,}  
\end{equation}
where $I$ is the $n\times n$ identity matrix.
The characteristic function $\chi$ of a linear autonomous DDE (of the
type considerd here) has at most finitely many roots with non-negative
real part \cite{HL93,S89}. This means that the stability theory of
equilibria of DDEs is very similar to that for ODEs.  An equilibrium
is stable if all eigenvalues of its linearisation have negative real
part. It changes its stability type when eigenvalues cross the
imaginary axis of the complex plane. When a single parameter is varied
there are two typical scenarios, which give rise to the standard local
bifurcations of equilibria of codimension one \cite{TS02,GH83,K04}:
\begin{itemize}
\item when a single real eigenvalue goes through zero the bifurcation
then one encounters, in the generic case, a \emph{saddle-node
bifurcation}, where two equilibria meet and disapper (or are created);
in the presence of additional properties of the equation (such as
symmetries) one encounters a \emph{branch bifurcation}, where two
additional equilibria bifurcate; \\[-6mm]
\item
when a complex conjugate pair of eigenvalues moves across the
imaginary axis of the complex plane then one encounters a \emph{Hopf
bifurcation} from which a small periodic orbit bifurcates.
\end{itemize}
These bifurcations of autonomous DDEs can be detected and tracked as
for ODEs by fixing the real parts of roots of the characteristic
function $\chi$ in \eqref{eq:chidef} to zero and freeing a second
parameter (such that $\eta\in\R^2$). In this way, one can compute
curves of saddle-node/branch and Hopf bifurcations in a
two-dimensional parameter space.  The curves divide the parameter
plane into regions where different numbers and types of equilibria
exist.

\subsection{Periodic orbits and their bifurcations}
\label{sec:perorbit}

\noindent
An initial function $x_0$ is a periodic point of the DDE
\eqref{eq:gendde} if it satisfies $x_{T_\Gamma}(s)=x_0(s)$ for all
$s\in[-\tau_2,0]$ and some $T_\Gamma\geq0$, where $T_\Gamma$ is called
the (minimal) period of $x_0$. Practically, one finds a periodic point
of \eqref{eq:gendde} as a \emph{periodic orbit}, which is a closed curve
$\Gamma:[0,1]\mapsto\R^n$ in the physical space satisfying the
\emph{periodic boundary-value problem}
\begin{equation}
  \label{eq:bvpgam}
  \begin{split}
    \dot \Gamma(t)&=T_\Gamma
    f(t,\Gamma(t),\Gamma(t-\tau_1/T_\Gamma)_{\mod[0,1]},
    \Gamma(t-\tau_2/T_\Gamma)_{\mod[0,1]},\eta)
    \mbox{\quad for $t\in(0,1]$,}\\
    \Gamma(0)&=\Gamma(1)\mbox{.}
  \end{split}
\end{equation}
The notation $\Gamma(t-\tau_j/T_\Gamma)_{\mod[0,1]}$ means that for
combinations of $t$, $\tau_j$ and $T_\Gamma$ for which
$t-\tau_j/T_\Gamma$ is outside of the interval $[0,1]$ we choose the
appropriate integer $\nu$ such that $t-\tau_j/T_\Gamma+\nu\in[0,1]$
and evaluate $\Gamma(t-\tau_j/T_\Gamma+\nu)$. From any solution
$\Gamma$ of \eqref{eq:bvpgam} one can find periodic points $x_0$ of
period $T_\Gamma$ through the relationship
\begin{equation}
  q(s)=\Gamma(s/T_\Gamma)_{\mod[0,1]}\mbox{\quad for $s\in[-\tau_2,0]$.}
  \label{eq:gamma2q}
\end{equation}
When the DDE is autonomous then the period $T_\Gamma$ in
\eqref{eq:bvpgam} is one of the unknowns that must be found together
with $\Gamma$. This can be done in the same way as for ODEs by the
adding one scalar equation, a so-called \emph{phase condition}
\cite{D07}, which determines $T_\Gamma$ and makes the solution
$\Gamma$ of \eqref{eq:bvpgam} locally unique.  Good starting points
for periodic orbits of autonomous DDEs are Hopf bifurcations where
small-ampitude periodic orbits branch off. Their period $T_\Gamma$ is
determined by the imaginary part of the pair of purely complex
conjugate eigenvalues at the Hopf bifurcation.

The situation is different when the DDE is periodically forced. Then
there are no equilibria, and the simplest invariant object is a
periodic orbit $\Gamma$ whose period equals the forcing period $T_f$,
that is, $T_\Gamma=T_f$ is known in \eqref{eq:bvpgam}. In fact, each
equilibrium $x_0(s) = \bar x_0, \ s\in[-\tau_2,0]$ of the autonomous
system for zero forcing amplitude gives rise to a trivial periodic
orbit of the form $\Gamma_0(t) \equiv \bar x_0$.  The trivial periodic
orbit $\Gamma_0$ can then be continued as a periodic orbit $\Gamma$ of
the periodically-forced DDE when the forcing is increased from
$0$. For the special case that the forcing is parametric, as in the
ENSO model \eqref{eq:tzip}, one finds that $\Gamma(t) = \Gamma_0(t)
\equiv 0$ remains constant for all forcing amplitudes; however, in
general, $\Gamma(t)$ is a non-constant periodic function when the
forcing is nonzero; note that generically the period of any periodic
orbit of a periodically-forced DDE is an integer multiple of the
forcing period.  Periodic orbits can generally not be found
analytically, but the boundary value
problem \eqref{eq:bvpgam} can be solved numerically with the
same numerical methods (for example, with Gauss collocation) that are
used to compute periodic orbits of ODEs \cite{RS07}.

The stability of a periodic orbit $\Gamma$ is determined by its
Floquet multipliers $\mu_i$, which are the eigenvalues of the
time-$T_\Gamma$ map of the point $x_0$ on the periodic orbit $\Gamma$
given by \eqref{eq:gamma2q}. If all $\mu_i$ are inside the unit circle
of the complex plane then $\Gamma$ is stable; conversely, the eigendirections
of the Floquet multipliers outside the unit circle correspond to
unstable directions.  Similarly to the case of equilibria, for a
time-periodic DDE one can construct a characteristic function
$\chi_p(\mu)$ such that the Floquet multipliers of $\Gamma$ are the
roots of $\chi_p$ \cite{SSH06,SS11}.  Generally, the characteristic
function $\chi_p$ and its roots are found numerically, after one has
computed $\Gamma$ and its period $T_\Gamma$ \cite{IS11}.  For a DDE
there are infinitely many roots of $\chi_p$ and, hence, Floquet
multipliers $\mu_i$, but they are also discrete and have the origin of
the complex plane as their only accumulation point; in particular,
there are at most finitely many unstable Floquet multipliers.  As for
ODEs, a periodic orbit changes its stability when Floquet multipliers
cross the unit circle, and this gives rise to the same well-known
local bifurcations of periodic orbits.  When a single parameter is
varied there are three typical scenarios, which give rise to the
standard generic local bifurcations of periodic orbits of
codimension-one \cite{TS02,GH83,K04}:
\begin{itemize}
\item
when a single real Floquet multiplier goes through $+1$ then one
encounters a \emph{saddle-node bifurcation of limit cycles}, where two
periodic orbits meet and disappear (or are created);\\[-6mm]
\item
when a single real Floquet multiplier goes through $-1$ then one
encounters a \emph{period-doubling bifurcation}, where a periodic
orbit of twice the period bifurcates;\\[-6mm]
\item
when a complex conjugate pair of eigenvalues goes through the unit
circle at $e^{\pm 2\pi \alpha i}$ then one encounters a
\emph{Neimark-Sacker or torus bifurcation}, where a smooth invariant
torus bifurcates (provided that low-order resonances are avoided, that
is, $\alpha \neq 1$, $1/2$, $1/3$, $2/3$,
$1/4$, $3/4$).
\end{itemize}
An example where bifurcations of the trivial periodic orbit can be
found analytically is the delayed Mathieu equation (also a
parametrically forced system) \cite{IS02}. More generally, however, 
these bifurcations of DDEs need to be detected and tracked numerically as
for ODEs, by appending the above conditions on the Floquet multipliers
to the continuation of the periodic orbit. In this way, one can
compute curves of saddle-node of limit cycle, period-doubling and
torus bifurcations of DDEs in a two-dimensional parameter space.

\subsection{Resonances and quasiperiodicity}
\label{sec:resonance}

\noindent
Resonant periodic orbits are special periodic orbits that correspond
to the locking of two different oscillations of the system.  They are
organised in DDEs in \emph{resonance tongues} and \emph{resonance
  surfaces} in the same way as in ODEs; see 
Figure~\ref{fig:bif2d} later and, for example, \cite{A83,MP94} for the
general theory for ODEs.  Along a torus bifurcation curve in a
parameter plane a Floquet multiplier $e^{\pm 2\pi \alpha i}$ moves
along the unit circle, that is, the quantity $\alpha$ changes. When
$\alpha = k/\ell \in \Q$, a two-parameter family of $k:\ell$
resonant periodic orbits (a so-called resonance surface) branches
off. The projection of this surface onto the parameter plane is a
resonance tongue, which is a region in the parameter plane where
periodic orbits exist that are $k:\ell$ locked to the original
periodic orbit (or to the forcing period in a forced DDE, such as the
parametrically forced ENSO model). Inside the resonance tongue there
are initially two periodic orbits, one of which is always
unstable. The root point of the resonance tongue is attached to the
bifurcation curve. Close to the root point the two periodic orbits
form $k:\ell$ torus knots.  In a cross section transverse to the flow
one finds a period-$\ell$ orbit; the return map to the cross section
acts on it by mapping each point to its $k$th neighbour in the
counter-clockwise direction. If the denominator $\ell$ is also greater
than $4$ both periodic orbits lie on an invariant torus (close to the
root point this is guaranteed by general theory \cite{K04}), one
periodic orbit is stable within the torus, one is unstable.

Similarly, one finds resonance tongues that are rooted at a line in
the parameter plane where the forcing is zero. At such a root point of
the $k:\ell$ resonance tongue the unforced system has a periodic orbit
$\Gamma$ such that the ratio of its period $T_\Gamma$ and the forcing
frequency $T_f$ is $\alpha = k/\ell \in \Q$.

By contrast, when $\alpha$ is irrational, that is, $\alpha \not\in
\Q$, then there is a smooth curve in the parameter plane starting at
the corresponding point on the torus bifurcation curve, or line of
zero forcing, along which one finds a torus that is densely filled by
a quasiperiodic trajectory.  The overall scenario is one of a
complicated interplay between locked and quasiperiodic dynamics on the
bifurcating torus. In particular, if one plots the rotation number
along any parameter path that intersects the resonance tongues one
finds the well-known devil's staircase consisting of infinitely many
intervals where the rotation number is constant (inside the
resonance tongues), as also observed in simulations of
ENSO models \cite{JNG94,JNG96,TSCJ94}.

A resonance tongue is bounded near its root point by saddle-node of
limit cycle bifurcation curves of the two orbits existing inside.  An
alternative method to the simulations as done in
\cite{JNG94,JNG96,TSCJ94} is to compute the $k:\ell$ resonance surface
by using the periodic BVP \eqref{eq:bvpgam} with period $T_\Gamma=\ell
T_f$ and two free parameters. In this way one obtains a surface
consisting of all locked periodic orbits inside the resonance tongue.
This compuation can be started at any point with rational $\alpha$ on
the torus bifurcation curve or the line where the system is unforced.
The continuation of the locked periodic orbits is numerically robust,
and the resonance tongue is obtained by projection onto the parameter
plane. This method has been introduced in \cite{MP94,SP07} for the
computation of resonance surfaces in maps and forced ODEs; Appendix A
explains how it can be adapted and implemented for a general DDE of
the form \eqref{eq:gendde}.

\subsection{Ready-to-use software for the bifurcation analysis of
  DDEs}
\label{sec:software}

\noindent
\texttt{DDE-Biftool} is a collection of \texttt{Matlab} functions that
perform each necessary task and that can be called by the user. For
example, there exist functions for calculating the eigenvalues of a
given equilibrium, the Floquet multipliers of a periodic orbit, for
converting a periodic BVP into a large (discretised) nonlinear system,
for solving nonlinear equations with Newton iteration, for adapting
the mesh on a given periodic solution, and for tracking solutions of
nonlinear systems if one has one more free parameter than equations.
Curve tracking of equilibria in a single parameter, of codimension-one
bifurcations of equilibria in two parameters, and of periodic orbits
of autonomous DDEs in a single parameter have been
implemented. Codimension-one bifurcations of periodic orbits can be
detected, but tracking them requires one to define a suitable extended
system manually.  Because \texttt{DDE-Biftool} permits one to add an
arbitrary number of conditions and parameters to the system, it is
possible to track solution surfaces instead of curves. We make use of
this feature for the computation of resonance tongues; see 
Appendix A.

The package \texttt{knut} (formerly known as \texttt{PDDE-Cont})
is a stand-alone alternative to \texttt{DDE-Biftool}. 
Its numerical methods are summarised in the review article by
Roose and Szalai \cite{RS07} and in \cite{SSH06}.

\section{Bifurcation analysis of the ENSO model}
\label{sec:bifs}

The main conclusion of Tziperman \emph{et al} \cite{TCZXB98} was that
in model \eqref{eq:tzip} the peak of the ENSO is locked to the end of
the year, which corresponds to the minimum of the excitation
$\kappa(t)$. That is, if $\kappa(t)=k_0+d_k\sin(t\pi/6)$ such that the
period of $\kappa$ is $12$ (time is in unit of months), then the
maximum of the thermocline anomaly, $h$, occurs at $t\approx
8\ldots9+12 \ell$ where $\ell$ is the denominator of the resonance. In
\cite{TCZXB98} the parameters were chosen such that $\ell=3$, so the
ENSO peak occured every three years. The authors also stated that in
all their simulations they found the locking to the peak to the end of
the year to be extremely robust with respect to all parameters.

We demonstrate here how the phenomenon of delay-induced resonances can
be studied systematically by means of a bifurcation analysis of
\eqref{eq:tzip}. More specifically, we determine all $k:\ell$
resonance tongues with $\ell \leq 11$ in the relevant region of the
$(d_k,k_0)$-plane of mean $k_0$ versus annual variation $d_k$ of the
amplification factor in \eqref{eq:kappadef}; they are chosen as the
bifurcation parameters since the interaction between ocean and
atmosphere is known only qualitatively as expressed by the function
$G$ shown in Figure~\ref{fig:gdef}. In contrast, confidence in the
other parameters of \eqref{eq:tzip} is higher: the delays $\tau_1$ and
$\tau_2$ are known with good accuracy, and the relation between
damping and feedback effects can be estimated directly from
observations (involving only oceanic quantities). These other
parameters are fixed at the values from \cite{TCZXB98}; see
Table~\ref{tab:parameters}.

\subsection{The autonomous regime for $d_k=0$}
\label{sec:autbif1d}

\begin{figure}[t!]
  \centering
  \includegraphics[width=1\textwidth]{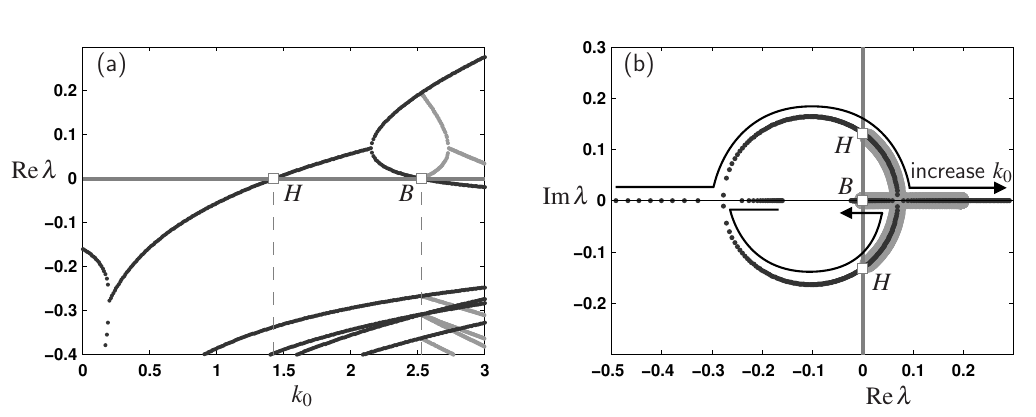}
  \caption{Stability of the trivial equilibrium $h \equiv 0$ and of
  bifurcating secondary equilibria of \eqref{eq:tzip} for $d_k=0$ as a
  function of $k_0$. Panel (a) shows the real parts of the eigenvalues
  of the equilibria as a function of $k_0$, and panel (b) shows the
  locus of the two leading eigenvalues in the complex plane; the
  eigenvalues of $h \equiv 0$ are plotted in black and those of the
  bifurcating equilibria in grey.}
  \label{fig:eigvals}
\end{figure}

\begin{figure}[t!]
  \centering
  \includegraphics[width=1\textwidth]{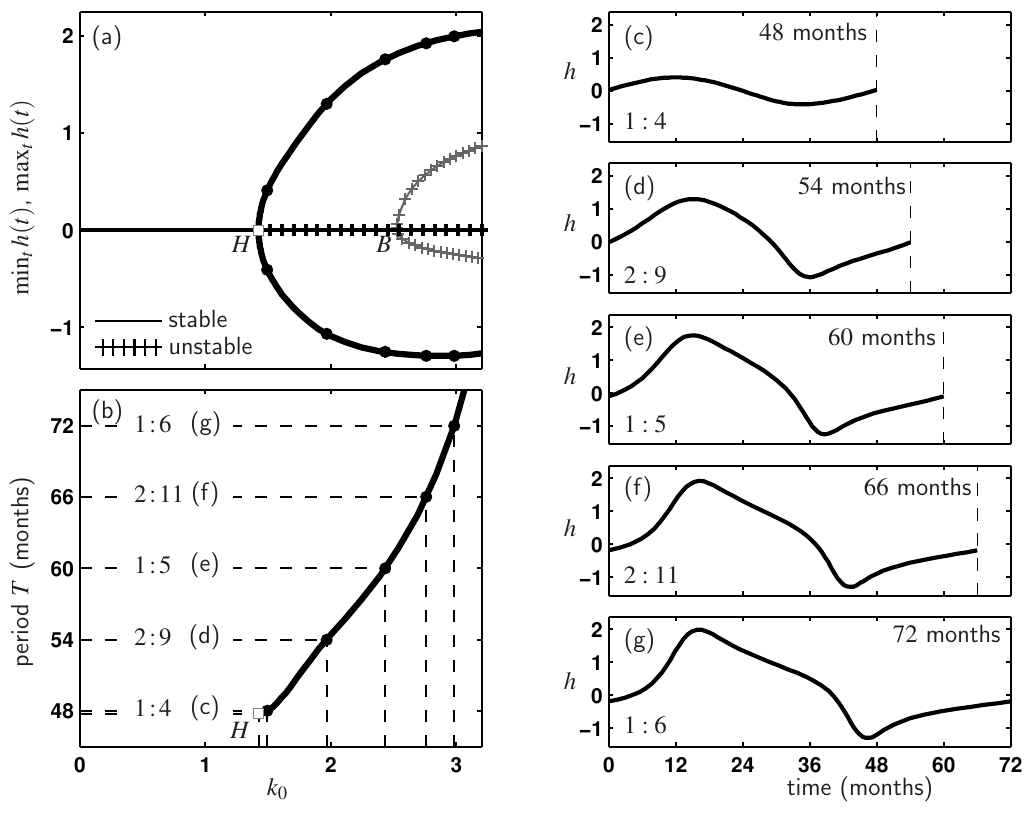}
  \caption{Single-parameter bifurcation diagram (a) of the autonomous
    system \eqref{eq:tzip} for $d_k=0$ and varying $k_0$. The
    equilibrium $h \equiv 0$ loses stability at the Hopf bifurcation
    $H$; a branch of secondary equilibria emerges from the branch
    point $B$. Panel (b) shows the period $T$ along the branch of stable
    periodic orbits that bifurcates from $H$, and panels (c)--(g) show
    time profiles of the indicated low-order resonant periodic orbits
    along the branch.}
  \label{fig:autbif1d}
\end{figure}

\noindent
The starting point of our bifurcation analysis is the autonomous case
where $d_k=0$. Then $h \equiv 0$ is an equilibrium irrespective of the
value of $k_0$, and we know that it is stable for $k_0=0$. We proceed
to continue this equilibrium with the package \texttt{DDE-Biftool}
\cite{ELR02} over the range $k_0 \in [0,3]$ while monitoring its
stability. The behaviour of the eigenvalues of $h \equiv 0$ is
illustrated by the black curves in Figure~\ref{fig:eigvals}.  Panel
(a) shows how the real part of the eigenvalues changes with $k_0$, and
panel (b) how the two leading eigenvalues move in the complex plane
with increasing $k_0$.  The trivial equilibrium $h \equiv 0$ is
initially stable and then loses its stability at $k_0 \approx 1.426$
in a Hopf bifurcation $H$, where a complex pair of eigenvalues crosses
the imaginary axis.  At $k_0 \approx 2.158$ the complex pair meets on
the real line of the complex plane and splits up into two real
eigenvalues. One of these real eigenvalues decreases and crosses the
origin of the complex plane at $k_0 \approx 2.526$ (the point is
labelled $B$ in Figures~\ref{fig:eigvals}
and~\ref{fig:autbif1d}). Physically, at $B$ the positive feedback
terms in \eqref{eq:tzip} exactly cancel the damping and the negative
feedback for $h=0$, and a branch of two nonzero equilibria bifurcates
and exists for $k_0$ past $B$. We continued these equilibria, and
their eigenvalues are shown as grey curves in
Figure~\ref{fig:eigvals}.  There is a pair of positive real
eigenvalues and, hence, the secondary equilibria are unstable in the
range $k_0 \in [2.526,3]$.

Figure~\ref{fig:autbif1d}(a) shows the one-parameter bifurcation
diagram for $k_0 \in [0,3]$. The equilibrium $h \equiv 0$ loses its
stability at $H$, and at the branch point $B$ the branch of secondary,
unstable equilibria emerges. Also shown in Figure~\ref{fig:eigvals}(a)
is the branch of bifurcating periodic orbits as computed by
continuation with \texttt{DDE-Biftool} from the Hopf bifurcation point
$H$. The periodic orbits are represented by plotting their minina and
maxima; they are all stable, as was checked with the computation of
the Floquet multipliers (not shown).  Of special interest are those
periodic orbits along the branch that give rise to low-order
resonances with the forcing period of $T_f=12$ months (once the annual
variation of the amplification is turned on, that is, $d_k$ is
increased from $0$).  Figure~\ref{fig:autbif1d}(b) shows the period
$T_\Gamma$ of the periodic orbits along the branch, where all $k:\ell$
resonant periodic orbits with $\ell \leq 11$ are indicated --- from
the $1:4$ resonance very close to the Hopf bifurcation $H$ to the
$1:6$ resonance at $k_0 \approx 2.982$. Panels (c)--(g) show the
respective time profiles of these resonant periodic orbits.

Figure~\ref{fig:autbif1d} gives a complete picture of the autonomous
dynamics of \eqref{eq:tzip} for $d_k=0$. Notice further from
Figure~\ref{fig:eigvals}(a) that all but the two leading eigenvalues
have quite large negative real parts. This suggests that the dynamics
of the DDE \eqref{eq:tzip} behaves like a two-dimensional ODE with
parametric periodic forcing for $k_0\leq3$ and delays that are
realistic for the ocean. For this reason, the bifurcation diagram
of the autonomous DDE is fairly simple. This is in contrast to the
much more complex bifurcation diagrams when the delays become large
enough to shift a large number of eigenvalues close to the imaginary
axis \cite{FSO07}.

\subsection{The dynamics with periodic forcing for $d_k>0$}
\label{sec:bif2d}

\begin{figure}[t!]
  \centering
  \includegraphics[width=1\textwidth]{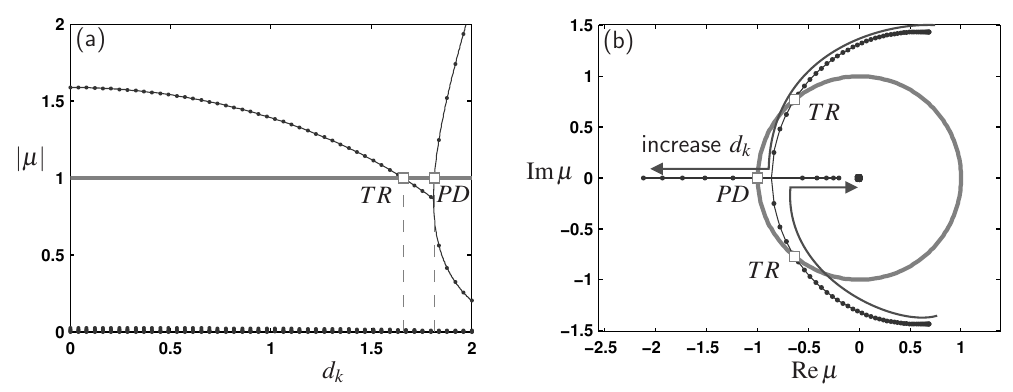}
  \caption{Stability of the trivial periodic solution $h \equiv 0$ as
  a function of the forcing $d_k$ for $k_0=1.8$.  Panel (a) shows the
  Floquet multipliers as a function of $d_k$, and panel (b) shows the
  locus of the two leading Floquet multipliers in the complex plane.}
  \label{fig:floqmul}
\end{figure}

\noindent
We now perform a bifurcation analysis of the periodically-forced DDE
\eqref{eq:tzip} for $d_k > 0$.  The trivial periodic solution $h
\equiv 0$ for $d_k=0$ gives rise for $d_k > 0$ to a periodic orbit
with period $T_f$.  Owing to the parametric nature of the periodic
forcing, this periodic orbit has zero amplitude, that is, it remains
the trivial solution $h \equiv 0$. As was the case for the equilibria
in the autonomous case, the main question is when $h \equiv 0$ is
stable. This can be determined by its continuation in $d_k$ for fixed
$k_0$. Figure~\ref{fig:floqmul} shows how the Floquet multipliers of
the trivial periodic orbit $h \equiv 0$ change as $d_k$ is increased
for fixed $k_0 = 1.8$, which is the value that was used in
\cite{TCZXB98}.  For $(k_0,d_k) = (1.8, 0)$ the periodic orbit $h
\equiv 0$ is initially unstable, since a single pair of complex
conjugate Floquet multipliers lies outside the unit circle. At $d_k
\approx 1.659$ this complex conjugate pair crosses the unit circle and
$h \equiv 0$ undergoes a torus bifurcation, labelled $TR$, and becomes
stable.  When $d_k$ is increased further, the complex conjugate pair
comes together on the real line of the complex plane at $d_k \approx
1.810$ and then splits up into two real Floquet multipliers. Very soon
thereafter, at $d_k \approx 1.814$, one of them decreases and crosses
the unit circle at $-1$ and, hence, $h \equiv 0$ undergoes a
period-doubling bifurcation, labelled $PD$, and becomes unstable again.

\begin{figure}[t!]
  \centering
  \includegraphics[width=1\textwidth]{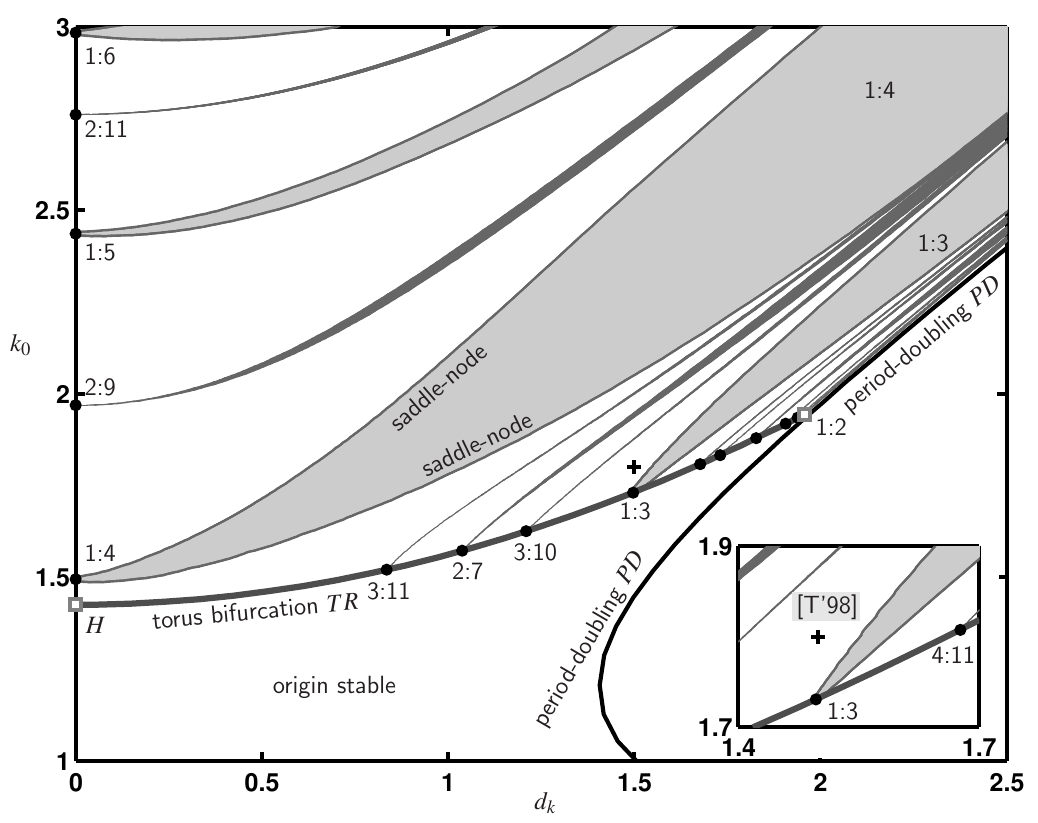}
  \caption{Two-parameter bifurcation diagram of the ENSO model
    \eqref{eq:tzip} in the $(d_k,k_0)$-plane $(d_k,k_0)$. All $k:\ell$
    resonance tongues with $\ell \leq 11$ are shown. They are rooted
    at points on the autonomous boundary $d_k=0$ and on the torus
    bifurcation curve $TR$, which ends on a period-doubling curve at a
    $1:2$ resonance. The trivial solution $h \equiv 0$ is stable in
    the region below the curve $TR$. The inset shows an enlargement
    near the tip of the $1:3$ resonance tongue; the point
    corresponding to the parameter values used by Tziperman \emph{et
    al} \cite{TCZXB98} is marked by a cross.}
  \label{fig:bif2d}
\end{figure}

We now consider the bifurcation diagram of \eqref{eq:tzip} in the
$(d_k,k_0)$-plane, which is shown in Figure~\ref{fig:bif2d}. The torus
bifurcation and the period-doubling bifurcation have been continued as
curves. The torus bifurcation curve $TR$ starts at the Hopf
bifurcation point on the automonous-limit line $d_k=0$ and it ends on
the period-doubling bifurcation curve $PD$ in a $1:2$ resonance. Below
the curve $TR$ and to the left of the curve $PD$ the trivial periodic
orbit $h \equiv 0$ is stable. To the right of the curve $PD$ and above
the curve $TR$ the periodic orbit $h \equiv 0$ is unstable.

\subsection{Resonance tongues in the ENSO model}
\label{sec:ensotongues}

\noindent
Above the torus bifurcation curve $TR$ one finds invariant tori and
associated $k:\ell$ resonance tongues, and Figure~\ref{fig:bif2d}
shows all resonance tongues for $\ell \leq 11$.  As we have already
seen in Figure~\ref{fig:autbif1d}, along the line $d_k=0$ and above
the Hopf bifurcation point $H$ we find five root points of these
resonance tongues. Additional root points of resonance tongues are
found along the torus bifurcation curve $TR$. Specifically, we observe
that the rotation number $\alpha$ gradually increases until the torus
bifurcation curve terminates in the $1:2$ resonance on the
period-doubling bifurcation curve $PD$; along $TR$ we find nine more
root points of $k:\ell$ resonance tongues with $\ell \leq 11$.  The
resonance tongues in Figure~\ref{fig:bif2d} have been computed as
resonance surfaces with the method described in Appendix A; they
were then projected onto the $(d_k,k_0)$-plane to yield the open
regions of locked periodic orbits and their boundaries.  Notice how, as theory
states \cite{A83,MP94}, the resonance tongues become very narrow with
increasing denominator $\ell$. We found that the corresponding $k:\ell$
locked periodic orbit of the ENSO model is stable and, hence,
observable practically everywhere throughout the considered parameter
region (exceptions are the parts of the narrow $3:11$, $3:10$,
$4:9$, and $5:11$ resonance tongues in the area where $d_k>2$).  The cross in
Figure~\ref{fig:bif2d} shows the parameter point chosen by Tziperman
\emph{et~al} \cite{TCZXB98}; it is slightly to the left of the $1:3$
tongue (see the inset), which explains that these authors observed
non-periodic motion closely resembling a period-three orbit.

\begin{figure}[t!]
  \centering
  \includegraphics[width=0.9\textwidth]{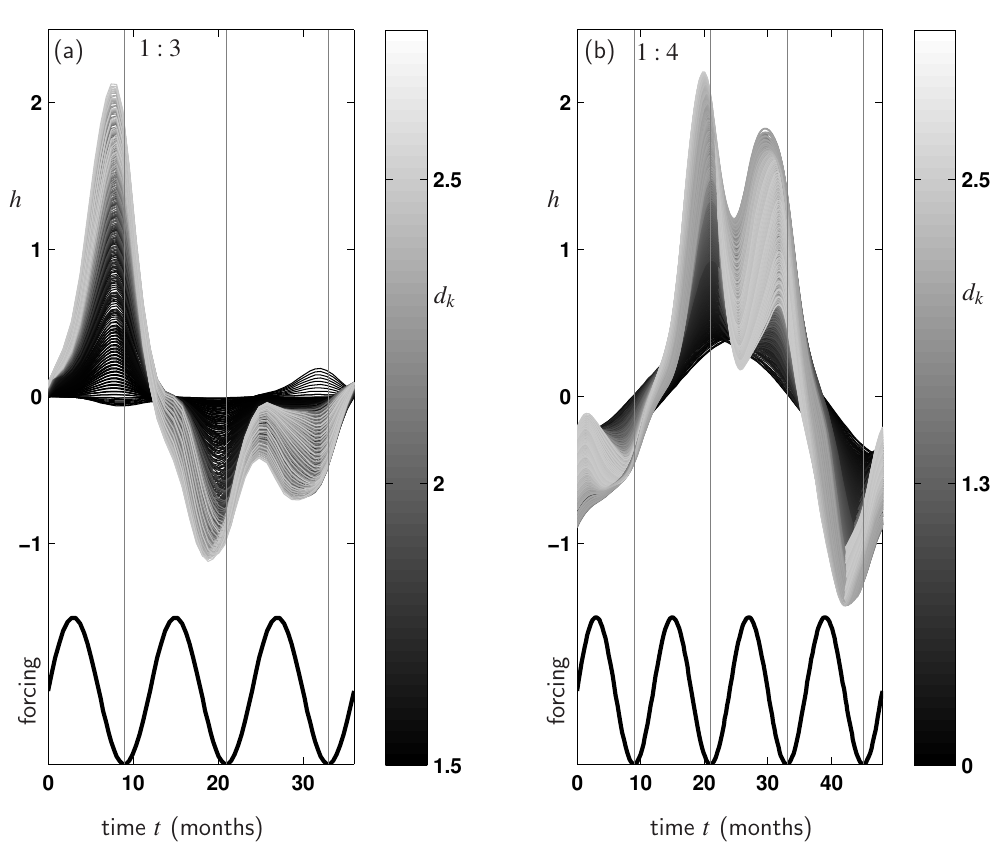}
  \caption{Locking of ENSO peak to the end of the year
    for all stable periodic orbits in the main tongues of $1:3$ (a) and
    $1:4$ resonance (b), where grayscale represents the 
    value of $d_k$ for which the periodic orbit is found. The black
    sinusoidal curve at the bottom indicates the phase of the forcing.}
  \label{fig:peaks}
\end{figure}

The other main feature observed 
in \cite{TCZXB98} is that the peak of the
oscillations appeared to be locked to the end of the year, that is,
the maximum of $h(t)$ occurs at time $t\approx 8\ldots9+12 \ell$.
To investigate how general this observation is 
Figure~\ref{fig:peaks} displays the profiles of all stable
periodic orbits in the main $1:3$ and $1:4$ resonance tongues,
together with the forcing. We find clear evidence of 
phase locking to the forcing in both cases, where the stable periodic
orbit reaches its maximum shortly before the forcing reaches its
minimum. Since the delays $\tau_1$ and $\tau_2$ enter
 directly into the phase of the 
forcing of the ENSO
model \eqref{eq:tzip},  it is clear that the two delays strongly influence 
the location of the maximum of $h$ relative to the forcing phase. 
Figure~\ref{fig:peaks} is evidence that the phase between $h$ and the
forcing  is 
nearly independent of the parameters $d_k$ and $k_0$. This 
observation provides a
certain justification for a reduction of the ENSO model even further
to contain only the delays and not the periodic forcing; see \cite{GZC08} 
for a demonstration of how resonances can be explored with Boolean
delay equations.

\section{Conclusions}
\label{sec:conc}

We presented a bifurcation study of how delay-induced resonances
organise the parameter plane of mean and variation of the annual
forcing of a model of ENSO with delayed feedback and forcing.  The
main message is that, colloquially speaking, the bifurcation theory
for the considered class of DDEs is exactly as that for ODEs, but
embedded into an infinite dimensional phase space.  Hence, equilibria,
periodic orbits and their bifurcations can be detected and tracked in
DDEs much as for ODEs, and software packages exist for this task. As
our case study demonstrates, numerical bifurcation anlysis is today
perfectly feasible for periodically-forced DDEs, and it allows one to
obtain insights beyond what is possible with long-time simulations.

Bifurcation analysis of equilibria and periodic orbits can give
evidence for the presence of deterministic chaos via the detection of
bifurcations associated with routes to chaos. A specific example is
chaotic dynamics due to the break-up of tori when resonance tongues
overlap.  In the ENSO model the resonance tongues do not overlap in
the parameter region we explored and, hence, no chaotic dynamics was
found.  However, when the forcing is large enough, resonance tongues
will generally overlap, and the method we used of computating
resonance tongues as surfaces is especially suitable in this
case. Hence, bifurcation analysis can add insight into the chaotic
dynamics in DDE models.

A limitation of bifurcation analysis is that chaotic behaviour can
only be studied indirectly, by finding bifurcations and instabilities
of equilibria and periodic orbits. In the study of climate systems
this imposes restrictions on the realism of the model. For example,
a realistic ocean boundary or ocean floor will introduce
fluctuations even at steady-state ocean flow. Thus, bifurcation
analysis in ocean and climate systems is typically performed on
idealised systems \cite{D08}. These idealized systems incorporate
the relevant physical effects but reduce realism such that steady
states become genuine equilibria. The Zebiak-Cane model for ENSO,
which provides the basis for the delayed oscillator model, is one
example of such idealised systems (described in the review
\cite{NBx98}).

Another issue in this regard is that --- even for
systems much simpler than climate models --- it is often very hard, if
not impossible, in practice to distinguish between deterministic chaos
and noise-induced fluctuations in observed data. One reason for this
may lie in a lack of hyperbolicity of the deterministic chaotic
attractor, which makes it non-robust (that is, sensitive to
arbitrarily small changes in the parameters); this non-robustness gets
`smeared out' if one adds random disturbances to the deterministic
model; see \cite{GCS08} for a study relevant to climate dynamics.
In particular,
for the ENSO these two hypotheses are still competing. Whereas
Tziperman \emph{et al} \cite{TCZXB98} presented the ENSO as if it
was mostly a deterministic chaotic trajectory, linear models subject
to noise are also used to explain, for example, the variation in
ENSO predictability \cite{TB01}. In terms of the two-parameter
bifurcation diagram in Figure~\ref{fig:bif2d}, linear models would
assume that the forcing parameters are from the stable region (below
curve $TR$), and that noise excites the system in the direction of the
dominant eigenvectors.

\section*{Acknowledgements}

The research of B.K. was supported by an FRDF grant from the Faculty
of Science, The University of Auckland, and that of J.S. by grant
EP/J010820/1 from the Engineering and Physical Sciences Research Council.


\vspace*{-1mm}

\appendix{Computation of resonance tongues as surfaces}
\label{sec:tongues}

\noindent
We now present technical details of how we compute a resonance surface
of a DDE of the form \eqref{eq:gendde}; our method is based on the
principle described in \cite{MP94,SP07} for ODEs and has been
implemented on top of \texttt{DDE-Biftool}.

\subsubsection{Branching off from a resonant torus bifurcation point}
\label{sec:trbranch}

We consider a $k:\ell$ resonant point $x$ on a previously computed torus
bifurcation curve in a parameter plane (such as the curve $TR$ in
Figure~\ref{fig:bif2d}), where the periodic orbit $\Gamma$
undergoing the torus bifurcation has a Floquet multiplier $e^{\pm 2\pi
\alpha i}$ with $\alpha = k/\ell$ and $0 < k < \ell$ are
coprime integers. We represent this point as the tuple
$(x,\alpha,\eta)$. Here $\eta$ is assumed to be of dimension two,
and the periodic orbit $\Gamma$ is the solution of the periodic BVP
\eqref{eq:bvpgam}.  Let $z(t)$ be the eigenfunction corresponding to
the Floquet multiplier $e^{2\pi \iu\alpha}$; then $z$ has period
$\ell$ (since $k$ and $\ell$ do not have a common divisor).  We now
choose a small radius $\rho$, and consider as initial guesses for the
initial circle (of radius $\rho$) of locked periodic orbits the family
\begin{displaymath}
  x_{0,\phi}(t)=\Gamma(\ell t)+
  \rho \cos \phi\re z(\ell t)+ \rho\sin\phi\im z(\ell
  t)\mbox{,\quad} 
  \phi\in[0,2\pi/\ell]\mbox{,\quad}t\in[0,1]\mbox{,}
\end{displaymath}
which is parametrised by the angle $\phi$. 
Note that here we have rescaled time again so that
$x_{0,\phi}$ has period $1$. 
We now perform a Newton iteration where we keep the two-dimensional
system parameter $\eta$ (initialised to its value at the resonance
point) free but fix the integral scalar products
\begin{equation}\label{eq:tanini}
    0=\int_0^1 \re z(\ell t)^T[x_\phi(t)-x_{0,\phi}(t)]\d t\mbox{,\quad}
    0=\int_0^1 \im z(\ell t)^T[x_\phi(t)-x_{0,\phi}(t)]\d t\mbox{.}
\end{equation}
Together with requiring that $x_\phi$ satisfies the periodic
boundary-value problem (similar to \eqref{eq:bvpgam} but with
$T_\Gamma=\ell T_f$)
\begin{equation}\label{eq:lockedbvp}
  \begin{split}
    \dot x_\phi(t)&=
    \ell T_f\,f(t,x_\phi(t),x_\phi(t-\tau_1/(\ell T_f)),
    x_\phi(t-\tau_2/(\ell T_f)),\eta)\mbox{,}\\  
    x_\phi(0)&=x_\phi(1)\mbox{,}
  \end{split}
\end{equation}
system \eqref{eq:tanini}--\eqref{eq:lockedbvp} is an equation for
the variables $(x_\phi,\eta)$ with a locally unique solution for every
$\phi$ and (sufficiently small) initial radius $\rho$.

\subsubsection{Branching off from an autonomous periodic orbit}
\label{sec:autbranch}

Creating the initial circle for a resonance surface starting from an
autonomous oscillation for zero forcing is slightly different from the
case of branching off from a torus curve. Suppose that at
$\eta=(\eta_1,\eta_2)=(\eta_1,0)$ the period $T_\Gamma$ of a periodic
orbit $\Gamma$ and the forcing period $T_f$ have a ratio $k/\ell$,
where $0 < k < \ell$ are again coprime; here we assume that the second
component $\eta_2$ of the parameter $\eta$ is the amplitude of the
forcing.  Then we choose as our initial guesses
$x_{0,\phi}(t)=x(\ell(t+\phi/(2\pi)))$ (again, rescaling time such
that the period of $x_{0,\phi}$ equals $1$). Two additional conditions
are given by fixing the initial forcing $\eta_2$ to the small radius
$\rho$, and by fixing the phase of the initial solution $x_\phi$:
\begin{equation}
  0=\int_0^1\dot x_{0,\phi}(t)^Tx_\phi(t)\d t
  \mbox{,\quad}\eta_2=\rho\mbox{.}\label{eq:auto:tanini}
\end{equation}
Solving \eqref{eq:lockedbvp} combined with \eqref{eq:auto:tanini}
results in an initial topological circle of locked orbits
$(x_\phi,(\eta_1,\eta_2))$, parametrised by $\phi\in[0,2\pi/\ell]$.

\subsubsection{Growing the surface circle by circle}
\label{sec:surfgrow}

Once we have an initial circle of solutions $x_\phi$ on the resonance
surface, we can continue the surface by computing nearby circles (which
are generally only topological circles). We parametrise the surface
locally by the angle $\phi$ on the circle and by its orthogonal
complement in the tangent space of the solution surface.
An equidistributed mesh of solutions on an already computed
topological circle is given. We pick a point $(x,\eta)_\mathrm{old}$
at angle $\phi$ on this circle. The tangent space in
$(x,\eta)_\mathrm{old}$ to the resonance surface is two-dimensional,
and can be computed as the two-dimensional nullspace of the
linearisation of \eqref{eq:lockedbvp} in $(x,\eta)_\mathrm{old}$. One
component of this tangent space, spanned by
$t_1=(x,\eta)_{\mathrm{tan},1}$, points along the old circle. We
denote its orthogonal complement in the tangent space by
$t_2=(x,\eta)_{\mathrm{tan},2}$, and take a small step (of size
$\delta$) from $(x,\eta)_\mathrm{old}$ along $t_2$ to obtain an
initial guess $(x,\eta)_\mathrm{new}$. Newton iteration is used to
solve \eqref{eq:lockedbvp}, augmented with the two linear equations
\begin{equation}
  \label{eq:surftan}
  0=\int_0^1 x_{\mathrm{tan},k}^T[x_\phi(t)-x_\mathrm{new}(t)]\d t
  +\eta_{\mathrm{tan},k}^T[\eta-\eta_\mathrm{new}]\mbox{\quad for $k=1$, $2$.}
\end{equation}
This gives for every $\phi\in[0,2\pi]$ a new solution $x_\phi$ on the
resonance surface. The family of new solutions $x_\phi$ forms again a
topological circle, so that the above procedure can be repeated to
`grow' the resonance surface as a family of topological circles. In
practice one has to perform computations only for angles
$\phi\in[0,2\pi/\ell]$, since for any $x_\phi$ the solution $t\mapsto
x_\phi(t-j/\ell)_{\mod[0,1]}$, $j=1\ldots \ell-1$, is also a solution.


\end{document}